\documentclass{elsart3-1}

\usepackage{amsmath,amssymb,amstext,dsfont,fancyvrb,float,fontenc,graphicx,subfigure, theorem}

\usepackage[english,francais]{babel}

\newtheorem{theorem}{Theorem}[section]
\newtheorem{lemma}[theorem]{Lemma}
\newtheorem{e-proposition}[theorem]{Proposition}

\newtheorem{e-definition}[theorem]{Definition\rm}

\newtheorem{theoreme}{Th\'eor\`eme}[section]

\newtheorem{proposition}[theoreme]{Proposition}

\newtheorem{definition}[theoreme]{D\'efinition\rm}

\renewcommand{\theequation}{\arabic{equation}}
\def\rmd{\mathrm{d}}
\def\rme{\mathrm{e}}
\def\rmi{\mathrm{i}}

\def\og{\leavevmode\raise.3ex\hbox{$\scriptscriptstyle\langle\!\langle$~}}
\def\fg{\leavevmode\raise.3ex\hbox{~$\!\scriptscriptstyle\,\rangle\!\rangle$}}

\journal{the Acad\'emie des sciences}
\begin{document}
\centerline{}
\begin{frontmatter}


\selectlanguage{english}
\title{Sample paths properties of Gaussian fields with equivalent spectral densities}


\selectlanguage{english}
\author[clausel]{Marianne Clausel}
\ead{marianne.clausel@imag.fr},
\author[vedel]{B\'eatrice Vedel}
\ead{vedel@univ-ubs.fr}


\address[clausel]{LJK,
Universit\'e de Grenoble--Alpes, CNRS
F38041 Grenoble Cedex 9}
\address[vedel]{LMAM, Universite de Bretagne Sud, Universit\'e Europ\'eene de Bretagne
Centre Yves Coppens
Bat. B, 1er et., Campus de Tohannic BP 573,
56017 Vannes, France.}


\medskip
\begin{center}
{\small Received *****; accepted after revision +++++\\
Presented by £££££}
\end{center}

\begin{abstract}
\selectlanguage{english}
We prove that if $X$ and $Y$ are two Gaussian fields with equivalent spectral densities, they have the same sample paths properties in any separable Banach space continuously embedded in $\mathcal{C}^0(K)$ where $K$ is a compact set of $\mathbb{R}^d$.

\vskip 0.5\baselineskip

\selectlanguage{francais}
\noindent{\bf R\'esum\'e} \vskip 0.5\baselineskip \noindent
{\bf Propri\'et\'es des trajectoires de champs gaussiens ayant des densit\'es spectrales \'equivalentes}
Nous montrons que si $X$ et $Y$ sont deux champs gaussiens \`a densit\'es spectrales \'equivalentes, ils ont m\^{e}me r\'egularit\'e dans tout espace de Banach s\'eparable s'injectant continument dans $\mathcal{C}^0(K)$ o\`u $K$ est un compact de $\mathbb{R}^d$.
\end{abstract}
\end{frontmatter}

\selectlanguage{english}


\section{Introduction}\label{s:intro}
In this note we are given two Gaussian random fields
$\{X(x)\}_{x\in\mathbb{R}^d}$ and $\{Y(x)\}_{x\in\mathbb{R}^d}$
both admitting stationary increments. We also assume that
these two fields admit a spectral density, that is there
exists two positive functions $f_X, f_Y\in
L^2(\mathbb{R}^d,(1\wedge |\xi|^2)\rmd \xi)$ such that,
\begin{equation}\label{e:defX}
X(x)=\int_{\mathbb{R}^{d}}(\rme^{\rmi x.\xi}-1)f_X^{1/2}(\xi)\rmd
\widehat{W}(\xi)\;,
\end{equation}
\begin{equation}\label{e:defY}
Y(x)=\int_{\mathbb{R}^{d}}(\rme^{\rmi x.\xi}-1)f_Y^{1/2}(\xi)\rmd
\widehat{W}(\xi)\;.
\end{equation}
We are also given $B$ a separable Banach space or a normed vector space, being the dual of a separable space. We assume that $B$ is continuously embedded in $\mathcal{C}^0(K)$ where $K$
denotes a compact of $\mathbb{R}^d$ (which is a separable Banach space). We aim at proving~:
\begin{theorem}\label{th:compreg}
Assume that there exists some $C>0$ such that
\begin{equation}\label{e:compdsp}
f_X(\xi)\leq C f_Y(\xi) \mbox{ for all }\xi\in\mathbb{R}^d\;.
\end{equation}
If the sample paths of $\{Y(x)\}_{x\in\mathbb{R}^d}$ a.s. belong
to $B$ then the sample paths of $\{X(x)\}_{x\in\mathbb{R}^d}$ a.s.
belong to $B$.
\end{theorem}
\section{Some classical results on probabilities in a Banach space}
Here $B_0$ is a separable Banach space. We denote $\mathcal{B}_0$ the Borel $\sigma-$ algebra  of $B_0$. If
$(\Omega, \mathcal{B}_{\Omega}, \mathbb{P})$ is a probability
space, a random element in $(B_0, \mathcal{B}_0)$ is a measurable
mapping from $(\Omega, \mathcal{B}_{\Omega}, \mathbb{P})$ in
$(B_0, \mathcal{B}_0)$.
\begin{definition}
Let $Z$ be a random element in $(B_0, \mathcal{B}_0)$ and $\mathbb{P}_Z$ its
distribution. A distribution of regular conditional
probability given $Z$ is a mapping $f\in B_0 \mapsto \mathbb{P}(\cdot|Z=f)$ such that :
\begin{enumerate}
\item $\forall f\in B_0$, $\mathbb{P}(\cdot|Z=f)$ is a probability measure on $\mathcal{B}$.\\
\item There exists a  $\mathbb{P}_{Z}$-negligible set $N$ such that
\[
\forall f\in B_0 \setminus N,\mathbb{P}(\Omega\setminus
Z^{-1}(f)|Z=f)=0\;.
\]
\item For all $A\in \mathcal{B}_{\Omega}$, the mapping
$f\mapsto \mathbb{P}(A|Z=f)$ is $\mathbb{P}_{Z}$-measurable and
\[
\mathbb{P}(A)=\int_{B_0}\mathbb{P}(A|Z=f)d\mathbb{P}_{Z}(f)\;.
\]
\end{enumerate}
\end{definition}
In separable Banach spaces, the distribution of
regular conditional probability given $Z$ exists and is unique.
More precisely~:
\begin{proposition}
For any random element $Z$ in
$(B_0, \mathcal{B}_0)$, there exists a distribution of conditional probability
given $Z$, $f \mapsto \mathbb{P}(\cdot|Z=f)$. If $f \mapsto \widetilde{\mathbb{P}}(\cdot|Z=f)$ is another one, then the set $\{f,\mathbb{P}(\cdot|Z=f)\neq\widetilde{\mathbb{P}}(\cdot|Z=f)\}$ is negligible.
\end{proposition}
\begin{definition}
A random element $X$ in $(B_0, \mathcal{B}_0)$ is Gaussian
if, for any linear form $L \in B_0^*$ (where $B_0^*$ denotes the dual
space of $B_0$), $L(X)$ is a real Gaussian random variable.
\end{definition}
The independence of Gaussian random elements is characterized as follows~\cite{27}~:
\begin{proposition}\label{Banachindep}
Two Gaussian random elements $X_1$ and $X_2$ in $(B_0, \mathcal{B}_0)$ are independent if for any linear forms $L_1$ and
$L_2$ of $B_0^*$, one has
$$
\mathbb{E}(L_{1}(X_{1})L_{2}(X_{2}))=0\;.
$$
\end{proposition}
\section{Proof of Theorem~\ref{th:compreg}}
The proof of Theorem~\ref{th:compreg} relies on the following
lemmas.
\begin{lemma}\label{lem:A1}
Let $X$ and $Y$ be two Gaussian fields of the form~(\ref{e:defX})
and~(\ref{e:defY}) with a.s. continuous sample paths. If $f_X \leq f_Y$ on $\mathbb{R}^{d}$, there exists two Gaussian
fields $X_1$ and $X_2$ with stationary increments, independent as
random elements with values in ${\mathcal{C}}^0(K)$, such that
$$
\{X(x)\}_{x\in\mathbb{R}^{d}}\overset{(\mathcal{L})}{=}\{X_{1}(x)\}_{x\in\mathbb{R}^{d}},\,
\{Y(x)\}_{x\in\mathbb{R}^{d}}\overset{(\mathcal{L})}{=}\{X_{1}(x)+X_{2}(x)\}_{x\in\mathbb{R}^{d}}\;.
$$
\end{lemma}

{\bf Proof.} Let us consider the Gaussian random field $Z$ defined
on $\mathbb{R}^d \times \mathbb{R}^2$ by its covariance function
\begin{eqnarray*}
&&\mathbb{E}(Z(x_1,\cdots,x_d;y_1,y_2)Z(x'_1,\cdots,x'_d;y'_1,y'_2))\\
&=&y_{1}y'_{1} \int_{\mathbb{R}^{d}}
(\rme^{\rmi x.\xi}-1)(\rme^{-\rmi x'.\xi}-1)f_{X}(\xi)\rmd\xi+y_{2}y'_{2}
\int_{\mathbb{R}^{d}} (\rme^{\rmi x.\xi}-
1)(\rme^{-\rmi x'.\xi}-1)(f_{Y}(\xi)-f_{X}(\xi))\rmd\xi\;.
\end{eqnarray*}
The inequality $f_Y-f_X\geq 0$ on $\mathbb{R}^d$ implies that
\[
\left((x_1,\cdots,x_d;y_1,y_2),(x'_1,\cdots,x'_d;y'_1,y'_2)\right)\mapsto\mathbb{E}(Z(x_1,\cdots,x_d;y_1,y_2)Z(x'_1,\cdots,x'_d;y'_1,y'_2))\;,
\]
is positive definite. Set now for any $x\in\mathbb{R}^{d}$, $X_{1}(x)=Z(x;1,0)$ and $X_{2}(x)=Z(x;0,1)$. Hence one has
\[
\{X(x)\}_{x\in\mathbb{R}^{d}}\overset{(\mathcal{L})}{=}\{X_{1}(x)\}_{x\in\mathbb{R}^{d}}
\mbox{ and }
\{Y(x)\}_{x\in\mathbb{R}^{d}}\overset{(\mathcal{L})}{=}\{Z(x,1,1)\}_{x\in\mathbb{R}^{d}}\overset{(\mathcal{L})}{=}\{X_{1}(x)+X_{2}(x)\}_{x\in\mathbb{R}^{d}}\;.
\]
Moreover, for all $x,x'$ in $\mathbb{R}^{d}$, $\mathbb{E}(X_{1}(x)X_{2}(x'))=0$. Using Proposition~\ref{Banachindep} and a Fubini theorem, since
the dual of $\mathcal{C}^{0}(B(0,1))$ is the set of Radon
measures, this last equality implies that $\{X_{1}(x)\}_{x\in\mathbb{R}^{d}}$ and
$\{X_{2}(x)\}_{x\in\mathbb{R}^{d}}$ are independent.

The assumptions and notations are now those of
Theorem~\ref{th:compreg}. The next lemma is a reformulation of
the Anderson inequality (see Theorem~11.9
of~\cite{22}) :
\begin{lemma}\label{lem:anderson}
Let $\{X(x)\}_{x\in K}$ a Gaussian random field defined on $K$
with a.s. continuous sample paths. Then, for any $r>0$ and $f \in
\mathcal{C}^{0}(K)$
$$
\mathbb{P}(\|X+f\|_{B}\leq r)\leq \mathbb{P}(\|X\|_{B}\leq r)\;.
$$
\end{lemma}
{\bf Proof.} Consider $X$ as a Gaussian random element in
$B_0=\mathcal{C}^0(K)$ which is a separable locally convex space.
Observe that in  Theorem~$9$ of~\cite{22} the set $C$ need only to be a convex, symmetric Borelian set of $B_0$ (personal communication of M. Lifshits). Hence, we can apply Theorem~$9$ of~\cite{22} to the Gaussian measure $\mathbb{P}_X$ and to the set $C=\{g\in B,\,\|g\|_{B}\leq r\}$ which is convex (since $\|\cdot\|_B$ is a norm), closed in $B$ since $B$ is either Banach either the dual of a Banach space and then a Borelian of $B_0=\mathcal{C}^{0}(K)$, and symmetric.

The following result can be deduced from
Lemma~\ref{lem:anderson}~:
\begin{lemma}\label{lem:A2}
Let $\{X_1(x)\}_{x\in K}$ and $\{X_2(x)\}_{x\in K}$ two
independent Gaussian random fields defined on a compact subset $K$
of $\mathbb{R}^{d}$ with a.s. continuous sample paths. For any $r>0$, one has
$$
\mathbb{P}(\| X_1+X_2 \|_B \le r) \le \mathbb{P}(\| X_1 \|_B \le
r)\;.
$$
\end{lemma}
{\bf Proof.} Since $X_1$ and $X_2$ are independent as random
elements in $\mathcal{C}^{0}(K)$, by definition of the conditional
probability, one has
\[
\mathbb{P}(\|X_1+X_2\|_B\leq r)= \int \mathbb{P}(\|X_1+f\|_B\leq r|X_2=f)\rmd\mathbb{P}_{X_2}(f)=\int \mathbb{P}(\|X_1+f\|_B\leq r)\rmd\mathbb{P}_{X_2}(f)\;.
\]
Lemma~\ref{lem:A1} applied to $X=X_1$
then implies that for any $f\in B$, $\mathbb{P}(\|X_1+f\|_B\leq r)\leq \mathbb{P}(\|X_1\|_B\leq r)$. Hence $\mathbb{P}(\|X_1+X_2\|_B\leq r)\leq \mathbb{P}(\|X_1\|_B\leq r)$.

Theorem~\ref{th:compreg} follows from these lemmas since the assumption
$f_X\leq Cf_Y$ and Lemma~\ref{lem:A1} imply
that
\[
\{Y(x)\}_{x\in\mathbb{R}^{d}}\overset{(\mathcal{L})}{=}\left\{\frac{1}{C^{1/2}}X(x)+X_{2}(x)\right\}_{x\in\mathbb{R}^{d}}\;.
\]
Lemma~\ref{lem:A2} then yields the required result.
{\footnotesize

\end{document}